\documentclass[11pt]{article}
\usepackage{amsfonts}
 \usepackage{bbm}
 \usepackage{amscd}
 \usepackage{xypic}
 \usepackage{amsmath}
 \usepackage{amssymb}
 \usepackage{amsthm}
 \usepackage{bbm}
 \usepackage{CJK}
 \usepackage{fancyhdr}
 \usepackage{hyperref}
 \usepackage{indentfirst}
 \usepackage{latexsym}
 \usepackage{mathrsfs}
 \allowdisplaybreaks
\usepackage{color}

\usepackage[top=1in,bottom=1in,left=1.25in,right=1.25in]{geometry}
 \def\<{\langle}
 \def\>{\rangle}
\def\a{\alpha}
\def\b{\beta}
\def\ci{\circ}
 \def\c{\cdot}
 
 \def\D{\Delta}

  \def\i{\iota}

\def\r{\rho}
\def\lr{\longrightarrow}

\def\o{\otimes}

 \def\v{\varepsilon}
 \def\vp{\varphi}
 
\def\<{\langle}
\def\>{\rangle}

\date{}
\begin{document}
\renewcommand{\baselinestretch}{1.2}
 \renewcommand{\arraystretch}{1.0}
\title{Characterization of Hopf Quasigroups}
\date{}
\author{{\it Wei WANG  and Shuanhong WANG  \footnote {Corresponding author:
 Shuanhong Wang, E-mail: shuanhwang@seu.edu.cn}}\\
{\small School of Mathematics, Southeast University, Jiangsu Nanjing 210096, China}\\
{\small E-mail:   weiwang2012spring@yahoo.com, shuanhwang@seu.edu.cn}}
\maketitle

 \vskip 0.5cm

 {\bf Abstract.} In this paper, we first discuss some properties of the Galois linear maps.
 We provide some equivalent conditions for  Hopf algebras
 and  Hopf  (co)quasigroups as its applications. Then let $H$ be a Hopf quasigroup with bijective antipode and
  $G$ be the set of all Hopf quasigroup automorphisms of $H$.
 We introduce a new  category $\mathscr{C}_{H}(\alpha,\beta)$ with $\alpha,\beta\in G$ over $H$
 and construct a new braided $\pi$-category $\mathscr{C}(H)$ with all the
 categories $\mathscr{C}_{H}(\alpha,\beta)$ as components.

 \vskip 0.5cm

{\bf Key words}:  Galois linear map; Antipode; Hopf (co)quasigroup; Braided $\pi$-category.
 \vskip 0.5cm

 {\bf Mathematics Subject Classification 2010:} 16T05.
\vskip 1.5cm

\section*{1. Introduction }

 The most well-known examples of Hopf algebras are the linear spans of (arbitrary)
 groups. Dually, also the vector space of linear functionals on a finite group carries the
 structure of a Hopf algebra. In the case of quasigroups (¡°nonassociative groups¡±)(see [A])
 however,  it is no longer a Hopf algebra but, more generally, a  Hopf quasigroup.
 \\

  In 2010,  Klim and Majid in [KM] introduced the notion of a Hopf quasigroup which
  was not only devoted to the development of this ¡¯missing link¡¯: Hopf quasigroup
   but also to understand the structure and relevant properties of the algebraic $7$-sphere,
   here Hopf quasigroups are not associative but the lack of this property is compensated by some axioms involving
  the antipode $S$.
  This mathematical object was studied in [FW] for twisted action as a generalization of Hopf
 algebras actions.
 \\

In 2000, Turaev  constructed a  braided monoidal category  in Freyd-Yetter categories of crossed
 group-sets (see [T1, T2]) which plays an important role of constructing some homotopy invariants.
 These crossed braided categories have been constructed in the setting of Hopf algebras (see [PS], [W], [ZD] and [LZ]).
 \\

Motivated by these two notions and structures, the aim of this paper is to construct classes of new
 braided $T$-categories from Hopf quasigroups.
 \\

Let $H$ be a Hopf quasigroup with bijective antipode. We denote by $G$ the set of all Hopf quasigroup automorphisms $\alpha$ of $H$ that
 satisfying $S\circ \alpha= \alpha \circ S$, and we consider $G$ a certain crossed product group $G\times G$.
\\

In Section 2, we recall some basic notations of algebras, coalgebras and bialgebras, we also introduce the definitions of Galois linear maps.
 \\

 In Section 3 and Section 4, we recall definitions and basic results related to Hopf (co) quasigroups,
 then we discuss some properties of the Galois linear maps in the setting of Van Daele [VD].
 As an application we provide some equivalent conditions for  Hopf algebras and  Hopf  (co)quasigroups.
 Then we also recall some notions related to Yetter-Drinfeld quasimodules over Hopf quasigroups and braided $\pi$-categories in Section 5 and Section 6.
\\

 In Section 5, we introduce a class of new categories $ \mathscr{C}_{H}(\alpha, \beta)$  of
 $(\alpha,\beta)$-Yetter-Drinfeld quasimodules associated with $\alpha,\beta\in G$ by generalizing  a Yetter-Drinfeld quasimodule over a Hopf quasigroup. Then in Section 6,  we prove $\mathscr{C}(H)$ is a monoidal
 cateogory and then construct a class of new braided $\pi$-categories $\mathscr{C}(H)$.
\\

\section*{2. Preliminaries }
\def\theequation{2. \arabic{equation}}
\setcounter{equation} {0} \hskip\parindent

An algebra $(A, m )$ is a vector space $A$ over a field $k$
 equipped with a map  $m : A\o A\lr A$. A unital  algebra $(A, m , \mu)$ is a vector space $A$ over a field $k$
 equipped with two maps $m : A\o A\lr A$ and $\mu : k\lr A$
 such that $m(id\o \mu ) =id=m(\mu \o id)$,
 where the natural identification $A\o k\cong k \cong k\o A$ is assumed. Generally, we write
 $1\in A$ for $\mu(1_k)$.
\\

 The algebra $(A, m , \mu)$ is
  called associative if $m(id \o m ) =m(m \o id)$. It is customary to write
$$
m(x\o y)=xy, \quad  \quad \quad \forall x, y\in C.
$$

A coalgebra $(C, \D )$ is a vector space $C$ over a field $k$
 equipped with a map  $\D : C\lr C\o C$. A counital  coalgebra $(C, \D , \v)$ is a vector space $C$ over a field $k$
 equipped with two maps $\D : C\lr C\o C$ and $\v : C\lr k$
 such that $(id\o \v )\D =id=(\v \o id)\D $,
 where the natural identification $C\o k\cong k \cong k\o C$ is assumed.

 The coalgebra $(C, \D , \v)$ is called coassociative if $(id \o \D )\D =(\D \o id)\D $. By using the Sweedler's
  notation [S], it is customary to write
$$
\D (x)=\sum x_{(1)}\o x_{(2)}, \quad   \quad \forall x\in C.
$$
\\

Given a counital coalgebra $(C, \D , \v )$ and a unital algebra $(A, m, \mu)$, the vector space $Hom(C, A)$ is a unital algebra with the
 product given by the convolution
\begin{eqnarray}
 (f*g)(x)=\sum f(x_{(1)})g(x_{(2)})
\end{eqnarray}
for all $x \in  C$, and unit element $\mu \v $. This algebra is denoted by $C*A$.
\\

Especially, we have the algebra $End(C)$ of endomorphisms on a given counital coalgebra $(C, \D , \v )$.
 Then we have the convolution algebra $C*End(C)$ with the unit element $Id: x\mapsto \v (x)id_C$.
 In case that the coalgebra $C$ is coassociative then $C*End(C)$ is an associative algebra.
 \\

A nonunital noncounital bialgebra  $(B, \D , m)$ is an algebra  $(B, m )$ and a coalgebra  $(B, \D )$
such that
$$
\D (xy)=\D (x)\D (y),\quad \quad \forall x, y\in B.
$$

A counital bialgebra  $(B, \D , \v, m )$ is a counital coalgebra  $(B, \D , \v)$ and an algebra $(B, m)$
such that
$$
\D (xy)=\D (x)\D (y), \quad \v (xy)=\v (x)\v (y), \quad \quad \forall x, y\in B.
$$

The multiplicative structure of a  counital bialgebra  $(B, \D , \v, m )$ is determined
 by the elements of $Hom(B, End(B))$:
 $$
  L: B\lr  End(B), a\mapsto L_a (L_a(x)=ax)
  $$
  and
  $$
  R :B \lr  End(B), a\mapsto R_a (R_a(x)=xa).
  $$
Obviously it suffices one of these maps to determine the multiplicative structure.
\\

A unital bialgebra  $(B, \D , m, \mu )$ is a coalgebra  $(B, \D )$ and a unital $(B, m, \mu)$
such that
$$
\D (xy)=\D (x)\D (y), \quad \D (1)=1, \quad \quad \forall x, y\in B.
$$

A unital counital  bialgebra  $(B, \D , \v,  m, \mu )$ is both a unital bialgebra  $(B, \D , m, \mu )$
 and a  counital bialgebra  $(B, \D , \v, m )$ such that $\v (1)=1$.
\\

Given a unital bialgebra  $(B, \D , m, \mu )$, we define the following two Galois linear maps
(cf.  [VD]):
\begin{eqnarray}
&& T_1: B\o B\lr B\o B, \quad \quad  T_1(x\o y)=\D (x)(1\o y);\\
&& T_2: B\o B\lr B\o B, \quad \quad  T_2(x\o y)=(x\o 1)\D (y)
\end{eqnarray}
for all $x, y \in  B$.
\\

It is easy to check $B\o B$ is  a left $B$-module and a right $B$-module
 with respective module structure:
$$
a(x\o y)=ax\o y, \quad \mbox {and } \quad (x\o y)a=x\o ya
$$
for all $a, x, y\in B$.
\\

Similarly, $B\o B$ is  a left $B$-comodule and a right $B$-comodule
 with respective comodule structure:
$$
\r _{H\o H}^l(x\o y)=\sum x_{(1)}\o (x_{(2)}\o y), \quad \mbox {and }
 \quad \r _{H\o H}^r(x\o y)=\sum (x\o y_{(1)})\o y_{(2)}
$$
for all $a, x, y\in B$.

\section*{3. Hopf algebras }
\def\theequation{3. \arabic{equation}}
\setcounter{equation} {0} \hskip\parindent

A  Hopf algebra $H$ is a unital associative  counital coassociative bialgebra  $(H, \D , \v,  m, \mu )$
 equipped with a linear map $S : H \lr H$ such that
\begin{eqnarray}
 \sum S(h_{(1)})h_{(2)}= \sum h_{(1)} S(h_{(2)})=\v(h)1
\end{eqnarray}
for all $h, g\in  H$.
\\

We have the main result of this section as follows:
\\

{\bf Theorem 3.1.}  Let $H:=(H, \D , \v,  m, \mu )$ be a unital associative  counital coassociative bialgebra.
  Then the following statements are equivalent:

 (1)   $H$  is a Hopf algebra;

 (2)  there is a  linear map $S : H \lr H$  such that $S$ and $id$ are invertible each other in
   the convolution algebra $H*H$;

 (3) the linear map $T_1: H\o H\lr H\o H$ is bijective, moreover,   $T_1^{-1}$ is a right $H$-module map
 and a left $H$-comodule map;

 (4) the linear map $ T_2: H\o H\lr H\o H$ is bijective, moreover,  $T_2^{-1}$ is a left $H$-module map
 and a right $H$-comodule map;

(5)  the element $L$ is invertible  in
   the convolution algebra $H * End (H)$;

(6)   the element $R$ is invertible  in
   the convolution algebra $H * End (H)$.
\\

{\bf Proof.} $(1)\Leftrightarrow (2)$. It follows from [S] that the part (1) is equivalent to the part (2).

  $(2)\Leftrightarrow (3)$. If the part (2) holds, it follows from [VD] that $T_1$  has the inverse $T^{-1}_1: A\o A\lr A\o A$
 defined by
 $$
 T^{-1}_1(a\o b)=\sum a_{(1)}\o S(a_{(2)})b
 $$
 for all $a, b\in H$.

It is not hard to check that $T_1^{-1}$ is a right $H$-module map
 and a left $H$-comodule map.

 Conversely, if the part (3) holds, then we introduce the notation, for all $a\in H$
 $$
 \sum a^{(1)}\o a^{(2)}:= T^{-1}_1(a\o 1).
 $$

Define a linear map $S: H\lr H$ by
$$
S(a)=(\v \o 1)\sum a^{(1)}\o a^{(2)}=\sum \v (a^{(1)}) a^{(2)}.
$$

Since $T_1^{-1}$ is a left $H$-comodule map, one has $(\r _{H\o H}^l\o id)T^{-1}_1=(id \o T^{-1}_1)\r ^l_{H\o H}$.
 That implies that, for all $a\in H$,
 $$
 \sum a^{(1)}{}_{(1)}\o a^{(1)}{}_{(2)}\o a^{(2)}=\sum a_{(1)}\o a_{(2)}{}^{(1)}\o a_{(2)}{}^{(2)}.
 $$

Applying $(id \o \v \o id )$ to the above equation, one obtains that
$$
 T^{-1}_1(a\o 1)=\sum a^{(1)}\o a^{(2)}=\sum a_{(1)}\o S(a_{(2)}).
 $$

Again, since $T_1^{-1}$ is the inverse of $T_1$  and it is a right $H$-module map,
 one concludes that
$$
a\o b=T^{-1}_1T_1(a\o b)=T^{-1}_1(\D (a)(1\o b))=\sum a_{(1)}\o S(a_{(2)})a_{(3)}b
$$
and
 $$
a\o b=T_1T^{-1}_1(a\o b)=T_1(\sum a_{(1)}\o S(a_{(2)})b)=\sum a_{(1)}\o a_{(2)}S(a_{(3)})b
$$

 Applying the counit to the first factor and taking $b=1$, we get equations (3.1).
 Thus, $S$ is the required antipode on $H$.

 $(2)\Leftrightarrow (4)$. Similarly, it follows from [VD] that
   $T_2$  has the inverse $T^{-1}_2: A\o A\lr A\o A$
   given by
   $$
   T^{-1}_2(a\o b)=aS(b_{(1)})\o b_{(2)}
   $$
   for all $a, b\in H$.
 Obviously,  $T_2^{-1}$ is a left $H$-module map
 and a right $H$-comodule map.

One  introduces the notation, for all $a\in H$
 $$
 \sum a^{[1]}\o a^{[2]}:= T^{-1}_2(1\o a).
 $$

Define  a linear map
 $S': H\lr H$ by
 $$
 S'(a)=(1\o \v )\sum a^{[1]}\o a^{[2]}
 =\sum a^{[1]}\v (a^{[2]}).
 $$

 By following the programm of arguments on $S$ we have
 that $S'$ satisfies equations (3.1).

 Furthermore, we now compute, for all $a\in H$,
 $$
  S'(a)=\sum S'(a_{(1)})\v (a_{(2)})=\sum \underline{S'(a_{(1)}) a_{(2)}}S(a_{(3)})=\sum \v(a_{(1)})S(a_{(2)})=S(a)
 $$

 Therefore, we have $S=S'$ and they are the required antipode on $H$.

 $(1)\Leftrightarrow (5)$. By hypothesis $B * End(B)$ is an associative algebra.
  The element $L$ is invertible in this algebra
  if and only if there exists $L': B \lr End(B)$ such that
$$
\sum L'(a_{(1)}) L(a_{(2)}) = \v (a) id = \sum L(a_{(1)})L'(a_{(2)})
$$

That implies that, for all $a, b\in H$,
$$
\sum L'(a_{(1)})(a_{(2)}b) = \v (a)b= \sum a_{(1)}L'(a_{(2)})(b),
$$
and in that case the inverse $L'$ is unique.

Defining $S: B\lr B$ by $S(a)=L'(a)(e)$ and taking $b=1$ and comparing this equation with (3.1)
 we obtain the desired result about the existence and uniqueness of $S$.

 Similarly for  $(1)\Leftrightarrow (6)$.

 This completes the proof. \hfill $\blacksquare$
\\

Let $G$ be a semigroup with the unit $e$. Then $(G, G)=\{(g, h)\mid g, h\in G\}$ is also
 a semigroup with the product:
 $$
 (x, y)(g, h)=(xg, yh)
 $$
for all $x, y, g, h\in G$.
\\

{\bf Corollary 3.2.}  Let $G$ be a semigroup with the unit $e$.
  Then the following statements are equivalent:

 (1)   $G$  is a group;

 (2)  there is a   map $S : G \lr G$  such that $S(g)g=e=gS(g)$ for all $g\in G$;

 (3) the map $T_1: (G, G)\lr (G, G), (g, h)\mapsto (g, gh)$ is bijective;

 (4) the map $T_2: (G, G)\lr (G, G), (g, h)\mapsto (gh, h)$  is bijective;

(5)  there is  a map $Q: G\lr End(G)$ such that the element
 $L: G\lr End(G), L(g)=L_g$, for all $g\in G$,
 satisfies $Q(g)(g) = e = g Q(g)(e)$

(6)  there is a   map $P: G\lr End(G)$ such that the element
 $R: G\lr End(G), R(g)=R_g$, for all $g\in G$,
 satisfies $P(g)(e) g = e = P(g)(g)$.

\section*{4. Hopf (co)quasigroups }
\def\theequation{4. \arabic{equation}}
\setcounter{equation} {0} \hskip\parindent

Recall from [A] that (an inverse property) quasigroup (or "IP loop") as a set $G$ with a product,
 identity $e$ and with the property that for each $ u \in G$ there is $u^{-1} \in G$ such that
$$
u^{-1}(uv) = v, \, \, \, \,  (vu)u^{-1} = v, \, \, \, \, \forall v \in G
$$

A quasigroup is flexible if $u(vu) = (uv)u$ for all $u, v \in G$
 and alternative if also $u(uv) = (uu)v$, $u(vv) =
(uv)v$ for all $u, v \in G$. It is called Moufang if $u(v(uw)) = ((uv)u)w$ for all
 $u, v, w \in G$.
\\

Recall from [KM, Definition 4.1 ] that  a {\sl Hopf quasigroup} is a unital algebra $H$ (possibly nonassociative)
 equipped with algebra homomorphisms $ \D : H \lr H \o  H$, $\v  : H \lr k$
  forming a coassociative coalgebra and a map $S : H \lr H$ such that
\begin{eqnarray}
&&  \sum S(h_{(1)})(h_{(2)}g)= \sum h_{(1)} (S(h_{(2)})g)=\v(h)g,\\
&& \sum (gS(h_{(1)}))h_{(2)}= \sum (gh_{(1)}) S(h_{(2)})=\v(h)g
\end{eqnarray}
for all $h, g\in  H$.  Furthermore, a Hopf
 quasigroup H is called {\sl flexible} if
$$
\sum h_{(1)}(gh_{(2)}) = \sum (h_{(1)}g)h_{(2)}, \quad  \forall h, g \in H,
$$
and {\sl Moufang} if
$$
\sum h_{(1)}(g(h_{(2)}f)) =\sum ((h_{(1)}g)h_{(2)})f \quad \forall  h, g, f \in H.
$$

Hence a Hopf quasigroup is a Hopf algebra iff its product is associative.
\\

Dually, we have
\\

 A {\sl Hopf coquasigroup} is a unital associative algebra
 $H$ equipped with counital algebra homomorphisms $ \D : H \lr H \o  H$, $\v  : H \lr k$
  and linear map  $S : H \lr H$ such that
\begin{eqnarray}
&& \sum S(h_{(1)})h_{(2)(1)}\o h_{(2)(2)}=1\o h=\sum h_{(1)} S(h_{(2)(1)})\o h_{(2)(2)},\\
&&\sum h_{(1)(1)}\o S(h_{(1)(2)})h_{(2)}=h\o 1=\sum h_{(1)(1)}\o h_{(1)(2)}S(h_{(2)})
\end{eqnarray}
for all $h \in  H$.  Furthermore, a Hopf
 coquasigroup H is called {\sl flexible} if
$$
\sum h_{(1)}h_{(2)(2)}\o h_{(2)(1)} = \sum h_{(1)(1)}h_{(2)}\o h_{(1)(2)}, \quad  \forall h \in H,
$$
and {\sl Moufang} if
$$
\sum h_{(1)}h_{(2)(2)(1)}\o h_{(2)(1)}\o h_{(2)(2)(2)}=\sum h_{(1)(1)(1)}h_{(1)(2)}\o h_{(1)(1)(2)}\o h_{(2)} \quad \forall  h \in H.
$$

Let  $(A, m , \mu)$ be a unital  algebra. Assume that $T: A\o A\lr A\o A$ is a map.
 Then we can define the following two coproduct maps:
\begin{eqnarray*}
&& \D ^r_T: A\lr A\o A, \quad a\mapsto T(a\o 1); \\
&& \D ^l_T: A\lr A\o A, \quad a\mapsto T(1\o a).
\end{eqnarray*}
\\

{\bf Definition 4.1.} With the above notation. We say that $T$
 is left (resp. right) compatible with $ \D ^r_T$, if $T(a\o b)= \D ^r_T(a)(1\o b)$
 (resp.  $T(a\o b)= (a\o 1)\D ^r_T(b)$),
 for all $a, b\in A$;

 Similarly, one says that $T$
 is left (resp. right) compatible with $ \D ^l_T$, if $T(a\o b)= \D ^l_T(a)(1\o b)$
 (resp.  $T(a\o b)= (a\o 1)\D ^l_T(b)$),
 for all $a, b\in A$.
\\

Dually, let  $(C, \D , \v)$ be a counital  coalgebra.
 Let   $T: A\o A\lr A\o A$ be a map.
 Then one  can define the following two product maps:
\begin{eqnarray*}
&& m ^r_T: A\o A\lr A, \quad a\o b\mapsto (1\o \v )T(a\o b); \\
&& m ^l_T:  A\o A\lr A, \quad a\o b\mapsto (\v \o 1 )T(a\o b).
\end{eqnarray*}

{\bf Definition 4.2.} With the above notation. We say that $T$
 is left (resp. right) compatible with $ m ^r_T$, if $T(a\o b)= (m ^r_T\o 1)(1\o \D )(a\o b)$
 (resp.  $T(a\o b)= (1\o m ^r_T)(\D \o 1)(a\o b)$),
 for all $a, b\in A$;

 Similarly, one says that $T$
 is left (resp. right) compatible with $ m ^l_T$, if $T(a\o b)= (m ^l_T\o 1)(1\o \D )(a\o b)$
 (resp.  $T(a\o b)= (1\o m ^l_T)(\D \o 1)(a\o b)$),
 for all $a, b\in A$.
\\

We now have the main result of this section as follows.
\\

{\bf Theorem 4.3.}  Let $H:=(H, \D , \v,  m, \mu )$ be a unital  counital coassociative  bialgebra.
  Then the following statements are equivalent:

 (1)   $H$  is a Hopf quasigroup;

 (2) the linear map $T_1, T_2: H\o H\lr H\o H$ is bijective, and   $T_1^{-1}$ is left compatible with
   $ \D ^r_{T_1^{-1}}$  and  right compatible with $ m ^l_{T_1^{-1}}$,
   at same time, the map $ T_2: H\o H\lr H\o H$ is bijective, moreover,  $T_2^{-1}$ is right compatible with
   $ \D ^l_{T_2^{-1}}$ and  left compatible with $ m ^r_{T_2^{-1}}$;

(3) the elements $L$ and $R$ is invertible  in  the convolution algebra $H * End (H)$.
\\

{\bf Proof.} $(1)\Leftrightarrow (2)$.
   If the part (2) holds, similar to Theorem 3.1, it is easy to check  that $T_1$  has the inverse $T^{-1}_1: A\o A\lr A\o A$
 defined by $ T^{-1}_1(a\o b)=\sum a_{(1)}\o S(a_{(2)})b $
 for all $a, b\in H$. Then we have
 \begin{eqnarray*}
&& \D ^r_{T_1^{-1}}: A\lr A\o A, \quad a\mapsto T^{-1}_1(a\o 1)=\sum a_{(1)}\o S(a_{(2)}); \\
&& m ^l_{T_1^{-1}}: A\o A\lr A, \quad a\o b\mapsto (\v \o 1)T^{-1}_1(a\o b)=S(a)b.
\end{eqnarray*}

It is not hard to check that $T_1^{-1}$ is left compatible with
   $ \D ^r_{T_1^{-1}}$  and  right compatible with $ m ^l_{T_1^{-1}}$.
\\

 Conversely, if the part (2) holds, then we define a linear map $S: H\lr H$ by
$$
S(a)=(\v \o 1)\D ^r_{T_1^{-1}}(a).
$$

Since $T_1^{-1}$ is a right compatible with $ m ^l_{T_1^{-1}}$ and it  is left compatible with
   $ \D ^r_{T_1^{-1}}$, one has, for all $a, b\in H$,
\begin{eqnarray*}
&&  T^{-1}_1(a\o b)\\
&=& (1\o m^l_{T^{-1}_1})(\D \o 1)(a\o b)\\
&=& \sum a_{(1)}\o (\v \o 1)T^{-1}_1(a_{(2)}\o b)\\
&=& \sum [a_{(1)}\o (\v \o 1)\D ^r_{T^{-1}_1}(a_{(2)})](1\o b)\\
&=&\sum a_{(1)}\o S(a_{(2)})b.
\end{eqnarray*}

Since $T_1^{-1}$ is the inverse of $T_1$, one hand, we conclude that
\begin{eqnarray*}
&&\sum a_{(1)}\o a_{(2)}[S(a_{(3)})b]\\
&=& T_1(\sum a_{(1)}\o S(a_{(2)})b)\\
&=& T_1T^{-1}_1(a\o b)\\
&=& a\o b\\
&=&T^{-1}_1T_1(a\o b)\\
&=&T^{-1}_1(\D (a)(1\o b))\\
&=&\sum a_{(1)}\o [S(a_{(2)})a_{(3)}]b.
\end{eqnarray*}
 Applying the counit to the first factor, one obtains equation (4.1).

We define another linear map $S': H\lr H$ by
$$
S'(a)=(1\o \v )\D ^l_{T_2^{-1}}(a), \quad \quad \forall a\in H.
$$

Similar to discussing with $S$, we can get equation (4.2).
By doing calculation, we have, for all $a\in H$
$$
S(a)=\sum S(a_{(1)})\v (a_{(2)})=\sum [S(a_{(1)})a_{(2)}]S'(a_{(2)})=\sum \v (a_{(1)})S'(a_{(2)})=S'(a).
$$

Thus, $H$ is a Hopf quasigroup.

 $(1)\Leftrightarrow (3)$. Similar to  $(1)\Leftrightarrow (5)$ and  $(1)\Leftrightarrow (6)$ in Theorem 3.1,
  it is not hard to finish the proof.

This completes the proof.  \hfill $\blacksquare$
\\

{\bf Corollary 4.4.}  Let $G$ be a nonempty with a product and with the unit $e$.
  Then the following statements are equivalent:

 (1)   $G$  is a quaigroup;

 (2)  there is a   map $S : G \lr G$  such that $S(g)(gh)=h=(hg)S(g)$ for all $g, v\in G$;

 (3) the map $T_1: (G, G)\lr (G, G), (g, h)\mapsto (g, gh)$ is bijective;

 (4) the map $T_2: (G, G)\lr (G, G), (g, h)\mapsto (gh, h)$  is bijective;

(5)  there is  a map $Q: G\lr End(G)$ such that the element
 $L: G\lr End(G), L(g)=L_g$, for all $g, h\in G$,
 satisfies $Q(g)(g)h = h = hg Q(g)(e)$

(6)  there is a   map $P: G\lr End(G)$ such that the element
 $R: G\lr End(G), R(g)=R_g$, for all $g\in G$,
 satisfies $P(g)(e)(g)h = h = P(g)(gh)$.
\\
\\

\section*{5. $(\alpha, \beta)$- Yetter-Drinfeld quasimodules over a Hopf quasigroup}
 \def\theequation{5. \arabic{equation}}
 \setcounter{equation} {0}
\hskip\parindent

 In this section, we will define the notion of a Yetter-Drinfeld quasimodule over
 a Hopf quasigroup that is twisted by two Hopf quasigroup automorphisms as well as the notion of a Hopf quasi-entwining structure
 and how to obtain such structure from automorphisms of Hopf quasigroups.
\\

Let $H$ be a Hopf quasigroup. In this section we simply write $\D (a)=a_{1}\o a_2$ for
$\D (a)=\sum a_{(1)}\o a_{(2)}$  for all $a\in H$.
\\

Recall from [AFGS] that  a left {\sl Hopf quasimodule} over a Hopf quasigroup $H$
 is a vector space $M$ with  a linear map $\cdot: H\otimes M\rightarrow
 M$ such that $ 1_{H}\cdot m = m$ and $
 a_{1}\cdot(S(a_{2})\cdot m)=S(a_{1})\cdot (a_{2}\cdot m)\,\,=\,\,\varepsilon(a)m$
 for all $a\in H$ and $m\in M$.
\\

A  left-right Yetter-Drinfeld quasimodule $M=(M,\cdot, \rho)$ over  a Hopf quasigroup $H$
 is a left $H$-quasimodule  $(M,\cdot)$ and  is a right $H$-comodule $(M,\rho)$ satisfying
 the following equation:
 $$
\left\{
  \begin{array}{l}
  (a_{2}\cdot m)_{(0)}\otimes (a_{2}\cdot m)_{(1)}a_{1}=a_{1}\cdot m_{(0)}\otimes a_{2}m_{(1)},\\
 m_{(0)}\otimes m_{(1)}(ab)=m_{(0)}\otimes (m_{(1)}a)b,\\
 m_{(0)}\otimes a(m_{(1)}b)=m_{(0)}\otimes (am_{(1)})b,
  \end{array}
\right.
$$
 for all $a, ab\in H$ and $m\in M$. Here we write  $\rho(m)=m_{(0)}\otimes m_{(1)}, \forall m\in M$.
\\

 Let $M$ and $N$ be two left-right Yetter-Drinfeld quasimodules over $H$.
  We call a linear morphism $f: M\rightarrow N$ a left-right Yetter-Drinfeld quasimodule
 morphism if
 $f$ is both a left $H$-quasimodule map and a right $H$-comodule map.
  We use $\mathscr{C}_{H}$ denote the category of left-right
 Yetter-Drinfeld quasimodules over $H$.
\\

In what follows, let  $H$ be a Hopf quasigroup with the bijective antipode $S$  and let $G$ denote the set of all automorphisms of a Hopf
quasigroup $H$. \\

{\bf Definition 5.1.} Let $\alpha, \beta\in  G$. A left-right {\sl $(\alpha, \beta)$-Yetter-Drinfeld quasimodule} over $H$ is a vector
space $M$, such that $M$ is a left $H$-quasimodule(with notation $h\otimes m\mapsto h\cdot m$) and a right $H$-comodule(with notation $M\rightarrow
M\otimes H$, $m\mapsto m_{(0)}\otimes m_{(1)}$) and with the following compatibility condition:
 \begin{equation}
\r (h\cdot m)=h_{21}\cdot m_{(0)}\otimes (\beta(h_{22})m_{(1)})\alpha(S^{-1}(h_{1})), \end{equation} for all $h\in H$ and $m\in M$. We denote by $
\mathscr{C}_{H}(\alpha,\beta)$ the category of left-right $(\alpha, \beta)$-Yetter-Drinfeld quasimodules, morphisms being both $H$-linear and
$H$-colinear maps. \\

 {\bf Remark.} Note that, $\alpha$ and $\beta$ are bijective,
 algebra morphisms, coalgebra morphisms and commute with $S$.
 \\

{\bf Proposition 5.2.} One has that Eq.(5.1) is equivalent to
 the following equations:
\begin{eqnarray} h_{1}\c m_{(0)}\o \beta(h_{2})m_{(1)}=(h_{2}\cdot m)_{(0)}\otimes (h_{2}\cdot m)_{(1)}\alpha(h_{1}). \end{eqnarray}

{\bf Proof.}  To prove this propostion, we need to use the property of antipode of a Hopf quasigroup. that is,
$(gS(h_{1}))h_{2}=(gh_{1})S(h_{2})=g\varepsilon(h)$, for all $g,h\in H$.

Eq.(5.1)$\Longrightarrow$ Eq.(5.2).

We first do calculation as follows: \begin{eqnarray*}
 &&(h_{2}\cdot m)_{(0)}\otimes (h_{2}\cdot m)_{(1)}\alpha(h_{1})\\
 &\stackrel{(5.1)}{=}& h_{22}\cdot m_{(0)}\otimes ((\beta(h_{23})m_{(1)})\alpha(S^{-1}(h_{21})))\alpha(h_{1})\\
 &=& h_{22}\cdot m_{(0)}\otimes ((\beta(h_{23})m_{(1)})\alpha(S^{-1}(h_{21})))\alpha(h_{1})\\
 &=& h_{2}\cdot m_{(0)}\otimes ((\beta(h_{3})m_{(1)})\varepsilon(h_{1})\\
 &=& h_{1}\cdot m_{(0)}\otimes \beta(h_{2})m_{(1)}.
\end{eqnarray*}

 For Eq.(5.2) $\Longrightarrow$ Eq.(5.1), we have
\begin{eqnarray*} &&h_{2}\cdot m_{(0)}\otimes (\beta(a_{3})m_{(1)})\alpha(S^{-1}(a_{1}))\\ &\stackrel{(5.2)}{=}& (h_{3}\cdot m)_{(0)}\otimes ((a_{3}\cdot
m)_{(1)}\alpha(a_{2}))\alpha(S^{-1}(a_{1}))\\ &=& (a_{2}\cdot m)_{(0)}\otimes (a_{2}\cdot m)_{(1)}\varepsilon(a_{1})\\ &=& (a\cdot m)_{(0)}\otimes (a\cdot
m)_{(1)}. \end{eqnarray*} This finishes the proof. \hfill $\blacksquare$ \\

{\bf Example 5.3.}  For $\beta\in G$, define $H_{\beta}=H$ as vector space over a field $k$, with regular right $H$-comdule structure and
left $H$-quasimodule structure given by $h\cdot h'=(\beta(h_{2})h')S^{-1}(h_{1})$, for all $h,h'\in H$. And more generally, if $\alpha,\beta\in G$,
define $H_{\alpha,\beta}=H$ with regular right $H$-comodule structure and left $H$-module structure given by $h\cdot
h'=(\beta(h_{2})h')\alpha(S^{-1}(h_{1}))$, for all $h,h'\in H$. If $H$ is flexible, then we get $H_{\beta}\in \mathscr{C}_{H}(id, \beta)$ and if we
add $H$ is $(\alpha,\beta)$-flexible, that is \begin{eqnarray*} \alpha(h_{1})(g\beta(h_{2}))=(\alpha(h_{1})g)\beta(h_{2}), \end{eqnarray*} for all $g,h\in
H$and $\alpha,\beta\in G$ . Then $H_{\alpha,\beta}\in  \mathscr{C}_{H}(\alpha,\beta)$. \\

 Let $\alpha,\beta\in G$. An $H$-bicomodule algebra $H(\alpha,\beta)$ is defined  as follows; $H(\alpha,\beta)=H$ as
 algebras, with comodule structures
 \begin{eqnarray*}
&& H(\alpha,\beta)\rightarrow H\otimes H(\alpha,\beta),\quad h\mapsto h_{[-1]}\otimes h_{[0]}=\alpha(h_{1})\otimes h_{2},\\ &&  H(\alpha,\beta)\rightarrow
H(\alpha,\beta)\otimes H,\quad h\mapsto h_{<0>}\otimes h_{<1>}=h_{1}\otimes \beta(h_{2}).
 \end{eqnarray*}
 Then we also consider the  Yetter-Drinfeld quasimodules like $_{H(\alpha,\beta)}\mathcal{YDQ}^{H}(H)$.
\\

{\bf Propostion 5.4.} $\mathscr{C}_{H}(\alpha,\beta)=_{H(\alpha,\beta)}\mathcal{YDQ}^{H}(H)$.

{\sl Proof.} Easy to check. \hfill $\blacksquare$
\\

\section*{6.  Main result }
 \def\theequation{4. \arabic{equation}}
 \setcounter{equation} {0}
 \hskip\parindent

 Recall from [FY] that a {\sl monoidal category}is a sixtuple $\mathcal{C}=(\mathcal{C},\mathbb{I},\otimes,a,l,r)$,
 where $\mathcal{C}$ is a category, $\otimes: \mathcal{C}\times\mathcal{C}\rightarrow\mathcal{C}$
 is a bifunctor called
 the {\sl tensor product}, $\mathbb{I}$ is an object of $\mathcal{C}$ called the {\sl unit}, and natural isomorphisms $a$
 ,$l$ and $r$ are called the {\sl associativity constraint}, the {\sl left unit constraint} and the {\sl right unit constraint},respectively,
 subject to the following two axioms:

 For all $U,V,X,Y\in \mathcal{C},$
 the {\sl associativity pentagon} is

 \centerline{
 \xymatrix{
 & ((U\otimes V)\otimes X)\otimes Y\ar[dl]^-{a_{U,V,X}\otimes id_{Y}}\ar[dr]_-{a_{U\otimes V, X,Y}} & \\
 (U\otimes (V\otimes X))\otimes Y\ar[d]^-{a_{U,V\otimes X,Y}} & & (U\otimes V)\otimes (X\otimes Y)\ar[d]_-{a_{U,V,X\otimes Y}}\\
 U\otimes ((V\otimes X)\otimes Y)\ar[rr]^-{id_{U}\otimes a_{V,X,Y}} & & U\otimes (V\otimes (X\otimes Y))
 }
 }

 and {\sl triangle diagram} is

 \centerline{
 \xymatrix{
 (U\otimes \mathbb{I})\otimes V\ar[rr]^-{a_{U,\mathbb{I},V}}\ar[dr]_-{r_{U}\otimes id_{V}} & & U\otimes (\mathbb{I}\otimes V)\ar[dl]^-{id_{U}\otimes l_{V}}\\
 & U\otimes V &
 }
 }
 A monoidal categoey $\mathcal{C}$ is {\sl strict} when all the above constraints are identities.
\\

  Let $\pi$ be a group,
  a category $\mathcal{C}$ over $\pi$
  is called a {\sl crossed $\pi$-category } (see [T1, T2]) if it is a monoidal category endow with two conditions:

  (1) A family of subcategories $\{\mathcal{C}_{\alpha }\}_{\alpha \in \pi}$ such that $\mathcal{C}$ is the disjoint union of this family
  and such that $U\otimes V\in \mathcal{C}_{\alpha \beta }$, for any $U\in \mathcal{C}_{\alpha }$,
 $V\in \mathcal{C}_{\beta }$, here $\mathcal{C}_{\alpha}$ is called the $\alpha$th component of $\mathcal{C}$;

 (2) Let $Aut(\mathcal{C})$ be the group of invertible strict tensor functors from $\mathcal{C}$ to itself, a group homomorphism
 $\varphi: \pi \rightarrow Aut(\mathcal{C})$, $\alpha \mapsto \varphi_{\alpha}$,called {\sl conjugation}, and assume $\varphi_{\alpha}(\varphi_{\beta})=\varphi_{\alpha\beta\alpha^{-1}}$,
 for all $\alpha,\beta\in \pi$.

 Furthermore, $\mathcal{C}$ is called strict when it is
 strict as a monoidal category.
\\
 Now given $\alpha\in \pi$ and assume $V$ is an object in $\mathcal{C}_{\alpha}$, now the functor $\varphi_{\alpha}$ will
 be denoted by $^{V}(\cdot)$, as in Turaev (cf. [T1],[T2]) or also $^{\alpha}(\cdot)$. We set $^{\overline{V}}(\cdot):=^{\alpha^{-1}}(\cdot)$.
 Then we have $^{V}id_{U}=id_{V_{U}}$ and $^{V}(g\circ f)=^{V}g\circ ^{V}f$.
 We then also have following two equations by the homomorphic condition of $\varphi_{\alpha}$,
 \begin{eqnarray*}
 ^{V\otimes W}(\cdot)=^{V}(^{W}(\cdot)), \quad ^{\mathbb{I}}(\cdot)=^{V}(^{\overline{V}}(\cdot))=^{\overline{V}}(^{V}(\cdot))=id_{\mathcal{C}}.
 \end{eqnarray*}
 Since the functor $^{V}(\cdot)$ is strict for all $V\in \mathcal{C}$, we have $^{V}(f\otimes g)=^{V}f\otimes ^{V}g$, for any morphisms $f,g$ in $\mathcal{C}$,
 in particularly, $^{V}\mathbb{I}=\mathbb{I}$. We call this the {\sl left index notation}.
\\

 A {\sl braiding} of a crossed category $\mathcal{C}$ is
 a family of isomorphisms
 \begin{eqnarray*}
 \{ c_{U,V}:U\otimes V\rightarrow {}^UV\otimes U\}
 \end{eqnarray*}
 where $U$ and $V$ are objects in $\mathcal{C}$,
 satisfying the following conditions:\\
 (i)~For any arrow $f\in \mathcal{C}_{\a }(U, U')$ and
 $g\in \mathcal{C}(V, V')$,
 $$
 (({}^{\a }g)\o f)\circ c _{U, V}=c _{U' V'}\circ (f\o g).
 $$
 (ii)~For all $ U, V, W\in \mathcal{C},$ we have
 $$
 c _{U\o V, W}=a_{{}^{U\o V}W, U, V}\circ (c _{U, {}^VW}\o
 id_V)\circ a^{-1}_{U, {}^VW, V}\circ (\i _U\o c _{V, W})
  \circ a_{U, V, W},
  $$
    $$c _{U, V\o W}=a^{-1}_{{}^UV, {}^UW, U}
 \circ (\i _{({}^UV)}\o c _{U, W})\circ a_{{}^UV, U, W}\ci
 (c _{U, V}\o \i_W)\circ a^{-1}_{U, V, W},$$
 where $a$ is the natural isomorphisms in the tensor category
 $\mathcal{C}$.\\
 (iii)~For all $ U, V\in \mathcal{C}$ and $\b\in \pi$,
 $$ \vp _{\b }(c_{U, V})=c_{\vp _{\b }(U), \vp_{\b }(V)}.
 $$

 A crossed category endowed with a braiding is called
 a {\sl braided $\pi$-category} (cf. [T1], [T2], [PS], [W]).
 \\

In this section, we will construct a class of new braided $\pi$-categories
 $\mathscr{C}(H)$ over any Hopf quasigroup $H$ with bijective antipode in the above Turaev's meaning.
 \\

 Denote $G^2=G
 \times G$
 a group with the following product:
$$
(\alpha,\beta)\ast (\gamma, \delta)=(\alpha\gamma, \delta\gamma^{-1}\beta\gamma),
$$
for all $\alpha,\beta, \gamma, \delta \in G$.
 The unit of this group is $(id,id)$ and $(\alpha,\beta)^{-1}=(\alpha^{-1},
 \alpha\beta^{-1}\alpha^{-1})$.\\

{\bf Proposition 6.1.} Let  $M\in \mathscr{C}_{H}(\alpha,\beta)$
 , $N\in \mathscr{C}_{H}(\gamma,\delta)$, with
 $\alpha,\beta,\gamma,\delta \in G$. Then we have

 (1) $M \otimes N \in \mathscr{C}_{H}((\alpha,\beta)\ast (\gamma, \delta))$
 with structures as follows:
\begin{eqnarray*}
&&h\c (m \otimes n)=\gamma (h_{1})\c m \otimes \gamma^{-1}\beta\gamma(h_{2})\c n,\\
&&m\otimes n \mapsto (m_{(0)}\otimes n_{(0)})\otimes
n_{(1)}m_{(1)}.
\end{eqnarray*}
for all $m\in M,n\in N$ and $h\in H.$

(2) Define ${}^{(\alpha,\beta)}N=N$ as
 vector space, with structures: for all $n\in N$ and $h\in H.$
$$
h\rhd n=\gamma^{-1}\beta\gamma\alpha^{-1}(h)\c n,
$$
$$
 n\mapsto n_{<0>}\otimes n_{<1>}=n_{(0)}\otimes \alpha\beta^{-1}(n_{(1)}).
$$
 Then $${}^{(\alpha,\beta)}N \in \mathscr{C}_{H}((\alpha,\beta)\ast (\gamma,\delta)\ast (\alpha,\beta)^{-1}),$$ where
$((\alpha,\beta)\ast (\gamma,\delta)\ast (\alpha,\beta)^{-1})=(\alpha\gamma\alpha^{-1}, \alpha\beta^{-1}\delta\gamma^{-1}\beta\gamma\alpha^{-1})$
as an element in $G^2$.
\\

{\bf Proof.} (1) Let $h,g\in H$ and $m\otimes n\in M\otimes N$. We can prove $1\cdot (m\otimes n)=m\otimes n$ and $h_{1}\cdot(S(h_{2})\cdot (m\otimes
n))=S(h_{1})\cdot(h_{2}\cdot (m\otimes n))=(m\otimes n)$, straightforwardly.

This shows that $M \otimes N$ is a left $H$-quasimodule, the right $H$-comodule condition is straightforward to check.

Next, we compute the compatibility condition as follows:
\begin{eqnarray*}
&& (h_{2}\cdot (m\otimes n))_{(0)}\otimes (h_{2}\cdot (m\otimes
n))_{(1)}\alpha\gamma(h_{1})\\
&=& (\gamma(h_{2})\cdot m\otimes \gamma^{-1}\beta\gamma(h_{3})\cdot n)_{(0)}\otimes (\gamma(h_{2})\cdot m\otimes
\gamma^{-1}\beta\gamma(h_{3})\cdot n)_{(1)}\alpha\gamma(h_{1})\\
&=& (\gamma(h_{2})\cdot m)_{(0)}\otimes (\gamma^{-1}\beta\gamma(h_{3})\cdot
n)_{(0)}\otimes (\gamma^{-1}\beta\gamma(h_{3})\cdot n)_{(1)}((\gamma(h_{2})\cdot m)_{(1)}\alpha\gamma(h_{1}))\\
&=& \gamma(h_{1})\cdot m_{(0)}\otimes
(\gamma^{-1}\beta\gamma(h_{3})\cdot n)_{(0)}\otimes ((\gamma^{-1}\beta\gamma(h_{3})\cdot n)_{(1)}\gamma\gamma^{-1}\beta\gamma(h_{2}))m_{(1)}\\
&=& \gamma(h_{1})\cdot m_{(0)}\otimes \gamma^{-1}\beta\gamma(h_{2})\cdot n_{(0)}\otimes \delta\gamma^{-1}\beta\gamma(h_{3})(n_{(1)}m_{(1)})\\
&=& h\cdot(m_{(0)}\otimes n_{(0)})\otimes \delta\gamma^{-1}\beta\gamma(h_{2})(n_{(1)}m_{(1)}).
\end{eqnarray*}
Thus $M \otimes N \in \mathscr{C}_{H}(\alpha\gamma, \delta\gamma^{-1}\beta\gamma)$.

(2) Obviously, the equations above define a  quasi-module and a comodule action of $N$. In what follows, we show the compatibility condition:
\begin{eqnarray*} &&(h\rhd n)_{<0>}\otimes(h\rhd n)_{<1>}\\ &=& (\gamma^{-1}\beta\gamma\alpha^{-1}(h)\c n)_{(0)}\otimes
\alpha\beta^{-1}((\gamma^{-1}\beta\gamma\alpha^{-1}(h)\c n)_{(1)})\\ &=& \gamma^{-1}\beta\gamma\alpha^{-1}(h)_{2}\cdot n_{(0)}\otimes
\alpha\beta^{-1}((\delta(\gamma^{-1}\beta\gamma\alpha^{-1}(h)_{3})n_{(1)}) \gamma(S^{-1}(\gamma^{-1}\beta\gamma\alpha^{-1}(h)_{1})))\\ &=& h_{2}\rhd
n_{(0)}\otimes (\alpha\beta^{-1}\delta\gamma^{-1}\beta\gamma\alpha^{-1}(h_{3})\alpha\beta^{-1}(n_{(1)}))\alpha\gamma\alpha^{-1}(S^{-1}(h_{1}))\\ &=&
h_{2}\rhd n_{(0)}\otimes (\alpha\beta^{-1}\delta\gamma^{-1}\beta\gamma\alpha^{-1}(h_{3})n_{<1>})\alpha\gamma\alpha^{-1}(S^{-1}(h_{1})) \end{eqnarray*} for
all $n\in N$ and $h\in H,$ that is ${}^{(\alpha,\beta)}N \in
\mathscr{C}_{H}(\alpha\gamma\alpha^{-1},\alpha\beta^{-1}\delta\gamma^{-1}\beta\gamma\alpha^{-1})$ .
\hfill $\blacksquare$
\\

{\bf Remark 6.2.} (1)  Note that, if $M \in \mathscr{C}_{H}(\alpha,\beta),\,\, N \in \mathscr{C}_{H}(\gamma,\delta)$ and $P\in \mathscr{C}_{H}(\mu, \nu)$,
then $(M\otimes N)\otimes P = M\otimes (N\otimes P)$ as objects in $\mathscr{C}_{H}(\alpha\gamma\mu,\nu\mu^{-1}\delta\gamma^{-1}\beta\gamma\mu).$

(2) Let $M \in \mathscr{C}_{H}(\alpha, \beta),
 \ \ N \in \mathscr{C}_{H}(\gamma,\delta),\,\, \mbox {and}\, (\mu,\nu)\in G^2$. Then by the
 above proposition, we have:
 $$
 {}^{(\alpha, \beta)\ast (\mu,\nu)}N={}^{(\alpha, \beta)}({}^{(\mu,\nu)}N),
$$ as objects in $\mathscr{C}_{H}(\alpha \mu\gamma \mu^{-1}\alpha^{-1}, \alpha\beta^{-1}\mu\nu^{-1}\delta\gamma^{-1}\nu\mu^{-1}\beta \mu\gamma
\mu^{-1}\alpha^{-1})$ and $$ {}^{(\mu,\nu)}(M\otimes N)= {}^{(\mu,\nu)}M \otimes {}^{(\mu,\nu)}N, $$ as objects in $\mathscr{C}_{H}(\mu\alpha \gamma
\mu^{-1},
 \mu\nu^{-1}\delta\gamma^{-1}\beta\alpha^{-1}\nu\alpha \gamma \mu^{-1})$.
\\

 Let $H$ be a Hopf quasigroup and
 $G^2=G
 \times G$.
 Define $\mathscr{C}(H)$ as the
  disjoint union of all $\mathscr{C}_{H}(\alpha,\beta)$
 with $(\alpha,\beta)\in G^2$. If we endow $\mathscr{C}(H)$
 with tensor product shown in Proposition 6.1 (1),
 then $\mathscr{C}(H)$ becomes
 a monoidal category.
\\

 Define a group homomorphism
 $\,\,\varphi: G\rightarrow Aut(\mathscr{C}(H)),\quad
 (\alpha, \beta) \,\,\mapsto \,\,\varphi(\alpha,\beta)\,\,$
 on components as follows:
\begin{eqnarray*} \varphi_{(\alpha,\beta)}: \mathscr{C}_{H}(\gamma,\delta)&\rightarrow& \mathscr{C}_{H}((\alpha,\beta)\ast
(\gamma,\delta)\ast (\alpha,\beta)^{-1}),\\ \quad \quad \quad \quad \quad
 \quad \varphi_{(\alpha,\beta)}(N)&=& {}^{(\alpha,\beta)} N,
\end{eqnarray*} and the functor $\varphi_{(\alpha,\beta)}$ acts as identity on morphisms.
\\

  We now prove the main result of  this article.
\\

{\bf Theorem 6.3.} Let $H$ be a Hopf quasigroup and
 $G^2=G\times G$. Then the category $\mathscr{C}(H)$ is a braided $\pi$-category over $G$.
\\

{\sl Proof.} We will finish the proof of this theorem with the following five steps:

 Step 1.  Let $M \in \mathscr{C}_{H}(\alpha,\beta)$ and
 $N \in \mathscr{C}_{H}(\gamma,\delta)$.  Take
  ${}^{M}N={}^{(\alpha,\beta)}N$ as explained in Proposition 6.1 (2).
  Define a map
  $$c_{M, N}: M \otimes N \rightarrow {}^{M}N \otimes M
  $$ by
$$
 c_{M,N}(m\otimes
 n)=n_{(0)}\otimes \beta^{-1}(n_{(1)})\cdot m,
$$
 for all $m\in M,n\in N.$

Then we prove that $c_{M,N}$ is an $H$-module map.
 Take  $h\cdot(m\otimes n)=\gamma(h_{1})\c m\otimes \gamma^{-1}\beta\gamma(h_{2})\c n$
 as explained in Proposition 6.1.
 \begin{eqnarray*}
&&c_{M,N}(h\cdot(m\otimes n))\\ &=& c_{M,N}(\gamma(h_1)\cdot m\otimes \gamma^{-1}\beta\gamma(h_2)\cdot n)\\ &=& (\gamma^{-1}\beta\gamma(h_{2})\cdot
n)_{(0)}\otimes \beta^{-1}((\gamma^{-1}\beta\gamma(h_{2})\cdot n)_{(1)})\cdot (\gamma(h_{1})\cdot m)\\ &=& (\gamma^{-1}\beta\gamma(h_{2})\cdot
n)_{(0)}\otimes (\beta^{-1}((\gamma^{-1}\beta\gamma(h_{2})\cdot n)_{(1)})\gamma(h_{1}))\cdot m\\ &=& \gamma^{-1}\beta\gamma(h_{3})\cdot n_{(0)}\otimes
((\beta^{-1}\delta\gamma^{-1}\beta\gamma(h_{4})\beta^{-1}(n_{(1)}))S^{-1}(\gamma(h_{2})))\gamma(h_{1})\cdot m\\ &=& \gamma^{-1}\beta\gamma(h_{2})\cdot
n_{(0)}\otimes (\beta^{-1}\delta\gamma^{-1}\beta\gamma(h_{3})\beta^{-1}(n_{(1)}))\varepsilon(h_{1})\cdot m\\ &=& \gamma^{-1}\beta\gamma(h_{1})\cdot
n_{(0)}\otimes (\beta^{-1}\delta\gamma^{-1}\beta\gamma(h_{2})\beta^{-1}(n_{(1)}))\cdot m, \end{eqnarray*} on the other side, we have \begin{eqnarray*}
  h\cdot c_{M,N}(m\otimes n) &=& h\cdot (n_{(0)}\otimes \beta^{-1}(n_{(1)})\cdot m)\\
  &=& \alpha(h_{1})\rhd n_{(0)} \otimes (\beta^{-1}\delta\gamma^{-1}\beta\gamma(h_{2})\beta^{-1}(n_{(1)}))\cdot m.
\end{eqnarray*}

Step 2.  Similarly we can check that $c_{M, N}$ is an $H$-comodule map.

Step 3.  Let $ P \in \mathscr{C}_{H}(\mu,\nu)$. Then we will check the equation:
$$
c_{M\otimes N,P}=(c_{M,^{N}P}\otimes id_{N})\circ (id_{M}\otimes c_{N,P}).
$$

 Using equations $^M(^N P)=^{M\otimes N}P$ and $^M(N\otimes P)=^M N\otimes ^M P$ to have \begin{eqnarray*}
&&(c_{M,^{N}P}\otimes id_{N})\circ (id_{M}\otimes c_{N,P})(m\otimes n\otimes p)\\ &=& c_{M,^{N}P}(m\otimes p_{(0)})\otimes \delta^{-1}(p_{(1)})\cdot n\\
&=& p_{(0)(0)}\otimes \beta^{-1}\gamma\delta^{-1}(p_{(0)(1)})\cdot m\otimes \delta^{-1}(p_{(1)})\cdot n\\ &=& p_{(0)}\otimes
\beta^{-1}\gamma\delta^{-1}(p_{(1)(1)})\cdot m\otimes \delta^{-1}(p_{(1)(2)})\cdot n\\ &=& p_{(0)}\otimes
\gamma^{-1}\beta^{-1}\gamma\delta^{-1}(p_{(1)})\cdot (m\otimes n)\\ &=& c_{M\otimes N,P}(m\otimes n\otimes p). \end{eqnarray*}

Step 4.  Similar we can check the equation:
$$
c_{M,N\otimes P}=(id_{^{M}N}\otimes c_{M,P})\circ(c_{M,N}\otimes id_{P})
$$

Step 5. We prove that  the map $c_{M,N}$ defined by $c_{M,N}(m\otimes n)=n_{(0)} \otimes \beta^{-1}(n_{(1)})\cdot m$ is bijective with
 inverse:
 $${c}_{M,N}^{-1}(n\otimes m)=\beta^{-1}(S(n_{(1)}))\cdot m \otimes n_{(0)}.$$

In fact,  we can prove $c_{M,N}\circ c_{M,N}^{-1}=id$. For all $m\in M, n\in
  N$, we have
\begin{eqnarray*}
&&c_{M,N}\circ c_{M,N}^{-1}(n \otimes m)\\ &=&c_{M,N}(\beta^{-1}(S(n_{(1)}))\cdot m\otimes n_{(0)})\\ &=&n_{(0)(0)}\otimes \beta^{-1}(n_{(0)(1)})\cdot
(beta^{-1}(S(n_{(1)}))\cdot m)\\ &=&n_{(0)(0)}\otimes \beta^{-1}(n_{(0)(1)}S(n_{(1)}))\cdot m\\ &=&n_{(0)}\otimes \beta^{-1}(n_{(1)(1)}S(n_{(1)(2)}))\cdot
m\\ &=&n_{(0)}\otimes \varepsilon(n_{(1)})m\\ &=& n\otimes m.
\end{eqnarray*}

The fact that $ c_{M,N}^{-1}\circ c_{M,N}=id$ is similar.

Therefore,  the category $\mathscr{C}(H)$ is a braided $\pi$-category over $G$.

 This completes the proof of the main result. \hfill $\blacksquare$
\\

\vskip 1.cm

{\bf Acknowledgements}
  The work was partially supported by the NNSFs of China (NO. 11371088 , NO.11571173 and NO.11601078),
  the Fundamental Research Funds for the Central Universities  (NO. CXLX12-0067) and the NSF of Jiangsu
 Province (No. BK20171348).

\vskip 1.cm

\bf References \rm
\bigskip

 [A] A.A. Albert, {\it Quasigroups I}, Trans. Amer. Math. Soc. 54 (1943), 507-519.

 [AFGS] J. N. Alonso $\acute{A}$lvarez, J. M. Fern$\acute{a}$ndez Vilaboa,   R. Gonz$\acute{a}$lez Rodr$\acute{1}$guez,
   C. Soneira Calvo, {\it Projections and Yetter-Drinfel'd modules over Hopf (co)quasigroups}, J Algebra, 443(2015), 153-199.

 [FW] X.L.Fang, S.H.Wang, {\it Twisted smash product for Hopf quasigroups} J. Southeast Univ. (English Ed.), 27(3)(2011), 343-346.

 [FY] P.J. Freyd, D.N.Yetter, {\it Braided compact closed categories with applications to low-dimensional topology},  Adv Math,  77(1989), 156-182.

 [KM] J. Klim, S. Majid, {\it Hopf quasigroups and the algebraic 7-sphere}, J Algebra, 323(2010), 3067-3110.

 [LZ] D.W. Lu, X.H. Zhang, {\it Hom-L-R-smash biproduct and the category of Hom-Yetter-Drinfel'd-Long bimodules}, Journal of Algebra and Its Applications, 17(7)(2018), 1850133.

 [PS] P. Panaite, M. D. Staic, {\it Generalized (anti) Yetter-Drinfel'd modules as components of a braided T-category},
  Israel J Math,  158(2007), 349-366.

 [S] M. E. Sweedler, {\it Hopf Algebras}, Benjamin, New York, 1969.

 [T1] V. G. Turaev, {\it Crossed group-categories},  Arab. J. Sci. Eng. Sect. C Theme Issues,  33(2C)(2008), 483-503.

 [T2] V. G. Turaev, {\it Quantum Invariants of Knots and $3$-Manifolds},  de Gruyter Stud  Math  de Gruyter, Berlin 18, 1994.

 [VD]  A. Van Daele, {\it Multiplier Hopf algebras}, Trans. Am. Math. Soc. 342(2) (1994), 917-932.

 [W]  S. H. Wang, {\it Turaev group coalgebras and twisted Drinfeld double}, Indiana Univ Math J,  58(3)(2009),   1395-1417.

 [ZD] X.H. Zhang, L.H. Dong, {\it Braided mixed datums and their applications on Hom-quantum groups}, Glasgow Mathematical Journal, 60(1)(2018), 231-251.

\end{document}